# The notion of $\psi$-weak dependence and its applications to bootstrapping time series[*]

**Paul Doukhan**

*Laboratoire de Statistique du CREST, Timbre J340*
*3, avenue Pierre Larousse, F-92240 Malakoff, and*
*SAMOS, Statistique Appliquée et Modélisation Stochastique,*
*Université Paris 1, Centre Pierre Mendès France,*
*90 rue de Tolbiac, F-75634 Paris Cedex 13, France*
*e-mail:* doukhan@ensae.fr

and

**Michael H. Neumann**

*Friedrich-Schiller-Universität Jena*
*Institut für Stochastik*
*Ernst-Abbe-Platz 2, D-07743 Jena, Germany*
*e-mail:* mneumann@mathematik.uni-jena.de

**Abstract:** We give an introduction to a notion of weak dependence which is more general than mixing and allows to treat for example processes driven by discrete innovations as they appear with time series bootstrap. As a typical example, we analyze autoregressive processes and their bootstrap analogues in detail and show how weak dependence can be easily derived from a contraction property of the process. Furthermore, we provide an overview of classes of processes possessing the property of weak dependence and describe important probabilistic results under such an assumption.

**AMS 2000 subject classifications:** Primary 60E15; secondary 62E99.
**Keywords and phrases:** Autoregressive processes, autoregressive bootstrap, mixing, weak dependence.

## 1. Mixing vs. weak dependence

For a long time mixing conditions have been the dominating type of conditions for imposing a restriction on the dependence between time series data. Such conditions were introduced at the end of the fifties by Rosenblatt (1956) and by the Saint Petersburg school, due to Ibragimov (1962). The notion was mainly used in a systematic way by statisticians since this notion fits quite well with nonparametric techniques; see for example Rosenblatt (1985) for details. The monograph edited by Eberlein and Taqqu (1986) describes the state of the art at that time, perhaps with the omission of some Russian and French authors. Doukhan and Portal (1987) obtained strong convergence results for the empirical process. Many examples of mixing processes are discussed in Doukhan (1994) and this monograph also includes previous results by Doukhan, León and Portal (see the citations included). Furthermore, Dehling and Philipp (2002) related several anterior work concerned with the strong invariance principle. However, the greatest advance in the theory occured after Rio (1993)'s work; his covariance inequality yields indeed sharp central limit results, see for this the monograph by Rio (2000) as well as Doukhan, Massart and Rio (1994). We would like to mention also contributions by Shao (1988) (who obtained sharp results under $\rho$-mixing), Shao and Yu (1996) (improved on the empirical CLTs) and Merlevède and Peligrad (2002) who tightened and precised several of the previous results involved with the strong coupling conditions. A summary of the state of the art is given in the very recent monograph by Bradley (2007); see also the nice review paper, Bradley (2005).

---

[*]This is an original survey paper

On the other hand, it turns out that certain classes of processes which are of interest in statistics are not mixing although a successive decline of the influence of past states takes place; see for example Dedecker *et al.* (2007). The simplest example of such a process is an AR(1)-process, $X_t = \theta X_{t-1} + \varepsilon_t$, where the innovations are independent and identically distributed with $\mathbb{P}(\varepsilon_t = 1) = \mathbb{P}(\varepsilon_t = -1) = 1/2$ and $0 < |\theta| \leq 1/2$; see also Rosenblatt (1980). It is clear that this process has a stationary distribution supported on $[-2, 2]$, and for a process in the stationary regime, it can be seen from the equality $X_t = \varepsilon_t + \theta\varepsilon_{t-1} + \cdots + \theta^{t-s-1}\varepsilon_{s+1} + \theta^{t-s}X_s$ that a past state $X_s$ can always be recovered from $X_t$. (Actually, since $|\varepsilon_t| > |\theta||\varepsilon_{t-1}| + \cdots + |\theta|^{t-s-1}|\varepsilon_{s+1}| + |\theta|^{t-s}|X_s|$ it follows that $X_t$ has always the same sign as $\varepsilon_t$ which means that we can recover $\varepsilon_t$ and therefore $X_{t-1}$ from $X_t$. Continuing in this way we can finally compute $X_s$.) This, however, excludes any of the commonly used mixing properties to hold. On the other hand, $X_s$ loses its impact on $X_t$ as $t \to \infty$. Another simple example of a non-mixing process is the Galton-Watson model with immigration; see Dedecker *et al.* (2007).

Besides the somehow artificial example above, there are many other processes of this type which are of great interest in statistics. For example, for bootstrapping a linear autoregressive process of finite order, it is most natural to estimate first the distribution of the innovations by the empirical distribution of the (possibly re-centered) residuals and to generate then a bootstrap process iteratively by drawing independent bootstrap innovations from this distribution. Now it turns out that commonly used techniques to prove mixing for autoregressive processes fail; because of the discreteness of the bootstrap innovations it is in general impossible to construct a coupling of two processes with different initial values.

Inspired by such problems, Doukhan and Louhichi (1999) and Bickel and Bühlmann (1999) introduced alternative notions called weak dependence and $\nu$-mixing, respectively. Actually, these two concepts are similar; instead of considering all bounded test functions (as for mixing) they proposed to weaken the condition by considering regular functions only which then includes many processes which are not mixing. Borovkova, Burton and Dehling (2001) is also an interesting advance in this area, those authors deal with functions of mixing sequences. A slightly simplified version of Doukhan and Louhichi (1999)'s definition is given here:

**Definition 1.1.** A process $(X_t)_{t\in\mathbb{Z}}$ is called *$\psi$-weakly dependent* if there exists a universal null sequence $(\epsilon(r))_{r\in\mathbb{N}}$ such that, for any $k$-tuple $(s_1,\ldots,s_k)$ and any $l$-tuple $(t_1,\ldots,t_l)$ with $s_1 \leq \cdots \leq s_k < s_k + r = t_1 \leq \cdots \leq t_l$ and arbitrary measurable functions $g : \mathbb{R}^k \to \mathbb{R}$, $h : \mathbb{R}^l \to \mathbb{R}$ with $\|g\|_\infty \leq 1$ and $\|h\|_\infty \leq 1$, the following inequality is fulfilled:

$$|\mathrm{cov}\,(g(X_{s_1},\ldots,X_{s_k}), h(X_{t_1},\ldots,X_{t_l}))| \leq \psi(k,l,\mathrm{Lip}\,g,\mathrm{Lip}\,h)\,\epsilon(r).$$

Here $\mathrm{Lip}\,h$ denotes the Lipschitz modulus of continuity of $h$, that is,

$$\mathrm{Lip}\,h = \sup_{x\neq y} \frac{|h(x)-h(y)|}{\|x-y\|_{l_1}},$$

where $\|(z_1,\ldots,z_l)\|_{l_1} = \sum_i |z_i|$, and $\psi : \mathbb{N}^2 \times \mathbb{R}_+^2 \to [0,\infty)$ is an appropriate function.

*Remark* 1.

(i) In Bickel and Bühlmann (1999), another type of weak dependence, called $\nu$-mixing, was introduced. Similarly to Definition 1.1, uniform covariance bounds over classes of functions with smooth averaged modulus of continuity are required. Usually, examples of processes obey both notions of weak dependence. We think that it is sometimes easier to verify a condition of weak dependence as in Definition 1.1 which we prefer for this reason.
(ii) Dedecker and Prieur (2005) introduced another related notion, called $\widetilde{\varphi}$-weak dependence, which is particularly adapted to expanding dynamical systems.
(iii) For the special case of causal Bernoulli shifts with i.i.d. innovations, that is $X_t = g(\varepsilon_t,\varepsilon_{t-1},\ldots)$, Wu (2005) introduced other measures of dependence which are somewhat connected to the coupling idea below. The notion of weak dependence considered here is more general and seems also to include processes which cannot be represented as Bernoulli shifts (*e.g.* associated processes perhaps cannot be written as Bernoulli shifts).

*Remark* 2. (*Some classes of weak dependence*)
Specific functions $\psi$ yield variants of weak dependence appropriate to describe various examples of models:

- $\kappa$-weak dependence for which $\psi(u,v,a,b) = uvab$; in this case we simply denote $\epsilon(r)$ as $\kappa(r)$.
- $\kappa'$ (causal) weak dependence for which $\psi(u,v,a,b) = vab$; in this case we denote $\epsilon(r)$ as $\kappa'(r)$. This is the causal counterpart of $\kappa$ coefficients which we recall only for completeness.

- $\eta$-weak dependence, $\psi(u,v,a,b) = ua+vb$; in this case we write $\epsilon(r) = \eta(r)$ for brevity.
- $\theta$-weak dependence is a causal dependence which refers to $\psi(u,v,a,b) = vb$; we write $\epsilon(r) = \theta(r)$ (Dedecker and Doukhan (2003)) for this causal counterpart of $\eta$ coefficients.
- $\lambda$-weak dependence $\psi(u,v,a,b) = uvab + ua + vb$; in this case we write $\epsilon(r) = \lambda(r)$ (Doukhan and Wintenberger (2007)). Besides the fact that it includes $\eta$- and $\kappa$-weak dependence, this new notion of $\lambda$-weak dependence is convenient, for example, for Bernoulli shifts with associated inputs.

It turns out that the notion of weak dependence is more general than mixing and allows to treat, for example, also Markovian processes driven by discrete innovations as they appear with time series bootstrap. In the next section we consider as an instructive example linear autoregressive processes of finite order and a corresponding bootstrap version thereof. We will demonstrate that the desired property of weak dependence readily follows from a contraction property which is typical for such models under standard conditions on the parameters. The approach described there is also applicable to proving weak dependence for many other classes of processes. Section 3 contains further examples of processes for which some sort of weak dependence has been proved. In Section 4 we give an overview of available tools under weak dependence. In particular, we provide a Donsker invariance principle and asymptotics for the empirical process. Furthermore, we also give Lindeberg Feller central limit theorems for triangular arrays. Finally, we provide probability and moment inequalities of Rosenthal and Bernstein type. Proofs of the new results in Section 2 are deferred to a final Section 5.

## 2. Autoregressive processes and their bootstrap analogues

In this section we intend to give a brief introduction to the basic ideas commonly used for verifying weak dependence. Most parts in this section are specialized to autoregressive processes of finite order and their bootstrap analogues.

We consider first a general real-valued stationary process $(X_t)_{t\in\mathbb{Z}}$. A simple and in many cases the most promising way of proving a property of weak dependence is via contraction arguments. For probability distributions $P$ and $Q$ on $(\mathbb{R}^d, \mathcal{B}^d)$ with finite mean, we define the metric

$$d(P,Q) = \inf_{(X,Y):\, X\sim P, Y\sim Q} \mathbb{E}\|X-Y\|_{l_1},$$

For $d=1$ and the $L_2$ instead of the $L_1$ distance, we obtain Mallows distance; see Mallows (1972). It is well known that such distances are suitable for metrizing weak convergence, that is, $d(P_n, P) \xrightarrow[n\to\infty]{} 0$ implies $P_n \Longrightarrow P$; see e.g. Bickel and Freedman (1981). Similar distances have also been used in the context of Markov processes to derive convergence of stationary distributions from convergence of the conditional distributions; see e.g. Dobrushin (1970) and Neumann and Paparoditis (2007). The following lemma shows that closeness of the conditional distributions in the above metric gives rise to estimates for covariances.

**Lemma 2.1.** *Suppose that $(X_t)_{t\in\mathbb{Z}}$ is a real-valued stationary process. Furthermore, let $s_1 \leq \cdots \leq s_k \leq t_1 \leq \cdots \leq t_l$ be arbitrary and let $g : \mathbb{R}^k \to \mathbb{R}$ and $h : \mathbb{R}^l \to \mathbb{R}$ be measurable functions. Let $(X'_t)_{t\in\mathbb{Z}}$ be any other version of the process, where $\mathbb{X}'_{s_k} = (X'_{s_k}, X'_{s_k-1}, \ldots)$ is independent of $\mathbb{X}_{s_k} = (X_{s_k}, X_{s_k-1}, \ldots)$. If $g$ is bounded and $\mathbb{E}|h(X_{t_1}, \ldots, X_{t_l})| < \infty$, then*

$$|\mathrm{cov}\,(g(X_{s_1},\ldots,X_{s_k}), h(X_{t_1},\ldots,X_{t_l}))|$$
$$\leq \ \|g\|_\infty \mathrm{Lip}\, h \ \mathbb{E}d\big(\mathbb{P}^{X_{t_1},\ldots,X_{t_l}|\mathbb{X}_{s_k}}, \mathbb{P}^{X'_{t_1},\ldots,X'_{t_l}|\mathbb{X}'_{s_k}}\big).$$

This lemma shows that a property of weak dependence follows from a convergence of the conditional distributions as the time gap to the lagged variables tends to infinity. The latter property can often be shown by appropriate coupling arguments. To get reasonably tight bounds for the covariances, one has to construct versions of the process, $(X_t)_{t\in\mathbb{Z}}$ and $(X'_t)_{t\in\mathbb{Z}}$ where the '$s_k$-histories' $\mathbb{X}_{s_k}$ and $\mathbb{X}'_{s_k}$ are independent but where the corresponding next process values $X_t$ and $X'_t$ are close; see the example below. Note that there is a close connection to the notion of $\tau$-dependence introduced by Dedecker and Prieur (2004). According to their Lemma 5, the infimum of $\mathbb{E}d(\mathbb{P}^{X_t|\mathbb{X}_s}, \mathbb{P}^{X'_t|\mathbb{X}'_s})$ is actually equal to their coefficient $\tau(t-s)$. Dedecker and Prieur (2004) used such coupling arguments to derive exponential inequalities and other interesting results, with applications to density estimation.

In the rest of this section we restrict our attention to a real-valued autoregressive process $(X_t)_{t\in\mathbb{Z}}$, which obeys the equation

$$X_t = \theta_1 X_{t-1} + \cdots + \theta_p X_{t-p} + \varepsilon_t, \quad t \in \mathbb{Z}. \tag{2.1}$$

The innovations $(\varepsilon_t)_{t\in\mathbb{Z}}$ are assumed to be independent and identically distributed with $\mathbb{E}\varepsilon_t = 0$. Furthermore, we make the standard assumption that the characteristic polynomial $\theta(z) = 1 - \theta_1 z - \cdots - \theta_p z^p$ has no zero in the unit circle. It is well known that there exists then a stationary solution to the model equation (2.1). We will assume that the process $(X_t)_{t\in\mathbb{Z}}$ is in the stationary regime. Then this process can be represented as a causal linear process,

$$X_t = \sum_{k=0}^{\infty} \alpha_k \varepsilon_{t-k}, \tag{2.2}$$

where $\alpha_k = \sum_{j\leq k} \sum_{k_1+\cdots+k_j=k} \theta_{k_1} \cdots \theta_{k_j}$. Denote by $\xi_1, \ldots, \xi_p$ the roots of the characteristic polynomial $\theta$ and let $\rho = \min\{|\xi_1|, \ldots, |\xi_p|\}$. Then, for any $\epsilon > 0$, there exists a $K_\epsilon < \infty$ such that, with $\rho_\epsilon = (1+\epsilon)/\rho$,

$$|\alpha_k| \leq K_\epsilon\, \rho_\epsilon^k \quad \forall k \in \mathbb{N}; \tag{2.3}$$

see e.g. Brockwell and Davis (1991, p. 85).

Convergence of $\mathbb{E}d(\mathbb{P}^{X_{t_1},...,X_{t_l}|\mathbb{X}_{s_k}}, \mathbb{P}^{X'_{t_1},...,X'_{t_l}|\mathbb{X}'_{s_k}})$ as $t_1 - s_k \to \infty$ can now be shown by a simple coupling argument. For this purpose, we consider a second (stationary) version of the autoregressive process, $(X'_t)_{t\in\mathbb{Z}}$, where $\mathbb{X}'_{s_k}$ is independent of $\mathbb{X}_{s_k}$. Note that $(X'_t)_{t\in\mathbb{Z}}$ can also be written as a linear process,

$$X'_t = \sum_{k=0}^{\infty} \alpha_k \varepsilon'_{t-k}.$$

Independence of $\mathbb{X}'_{s_k}$ and $\mathbb{X}_{s_k}$ is equivalent to the fact that $\varepsilon_{s_k}, \varepsilon_{s_k-1}, \ldots$ and $\varepsilon'_{s_k}, \varepsilon'_{s_k-1}, \ldots$ are independent. On the other hand, we have some freedom to couple the innovations after time $s_k$. Here we only have to take care that both sequences $(\varepsilon_t)_{t\in\mathbb{Z}}$ and $(\varepsilon'_t)_{t\in\mathbb{Z}}$ consist of independent random variables. A reasonably good coupling is obtained by feeding both processes after time $s_k$ with one and the same sequence of innovations, that is, $\varepsilon'_{s_k+1} = \varepsilon_{s_k+1}, \varepsilon'_{s_k+2} = \varepsilon_{s_k+2}, \ldots$. This gives, for $t \geq s_k$, that

$$X_t - X'_t = \sum_{l=0}^{\infty} \alpha_{t-s_k+l}(\varepsilon_{s_k-l} - \varepsilon'_{s_k-l}). \tag{2.4}$$

Since $d(\mathbb{P}^{X_{t_1},...,X_{t_l}|\mathbb{X}_{s_k}}, \mathbb{P}^{X'_{t_1},...,X'_{t_l}|\mathbb{X}'_{s_k}}) \leq \mathbb{E}\bigl(\sum_{j=1}^{l} |X_{t_j} - X'_{t_j}| \bigm| \mathbb{X}_{s_k}, \mathbb{X}'_{s_k}\bigr)$ we obtain, in conjunction with (2.3), the following assertion:

**Lemma 2.2.** *Let $(X_t)_{t\in\mathbb{Z}}$ and $(X'_t)_{t\in\mathbb{Z}}$ be two versions of the autoregressive process as described above. If $\mathbb{E}|\varepsilon_0| < \infty$, then, for $s_k \leq t_1 \leq \ldots \leq t_l$,*

$$\mathbb{E}d\bigl(\mathbb{P}^{X_{t_1},...,X_{t_l}|\mathbb{X}_{s_k}}, \mathbb{P}^{X'_{t_1},...,X'_{t_l}|\mathbb{X}'_{s_k}}\bigr) \leq 2\, l\, K_\epsilon\, \frac{1}{1-\rho_\epsilon}\, \rho_\epsilon^{t_1-s_k}\, \mathbb{E}|\varepsilon_0|.$$

Lemma 2.1 and Lemma 2.2 imply the following weak dependence property.

**Corollary 2.1.** *Suppose that $(X_t)_{t\in\mathbb{Z}}$ is a stationary process satisfying the above conditions. Furthermore, let $s_1 \leq \cdots \leq s_k \leq t_1 \leq \cdots \leq t_l$ be arbitrary and let $g: \mathbb{R}^k \longrightarrow \mathbb{R}$ and $h: \mathbb{R}^l \longrightarrow \mathbb{R}$ be measurable functions. If $g$ is bounded and $\mathbb{E}|h(X_{t_1},\ldots,X_{t_l})| < \infty$, then*

$$|\mathrm{cov}\,(g(X_{s_1},\ldots,X_{s_k}), h(X_{t_1},\ldots,X_{t_l}))| \leq \|g\|_\infty \, \mathrm{Lip}\, h\, K_1\, l\, \rho_\epsilon^{t_1-s_k},$$

*where $K_1 = 2\,\frac{1}{1-\rho_\epsilon}\, K_\epsilon\, \mathbb{E}|\varepsilon_0|$.*

Now we define the autoregressive bootstrap. We assume that observations $X_{1-p}, \ldots, X_n$ are available. Let $\widehat{\theta}_n = (\widehat{\theta}_{n,1},\ldots,\widehat{\theta}_{n,p})'$ be any consistent estimator of $\theta = (\theta_1,\ldots,\theta_p)'$, that is, $\widehat{\theta}_n \xrightarrow{P} \theta$, as $n \to \infty$. (The least squares and the Yule-Walker estimator are even $\sqrt{n}$-consistent.) Let $\boldsymbol{X}_t = (X_{t-1},\ldots,X_{t-p})'$ be the vector of the $p$ lagged observations at time $t$. We define residuals

$$\widetilde{\varepsilon}_t = X_t - \boldsymbol{X}'_t \widehat{\theta}_n$$

and re-center them as

$$\widehat{\varepsilon}_t = \widetilde{\varepsilon}_t - \frac{1}{n}\sum_{t=1}^{n} \widetilde{\varepsilon}_t.$$

Now we draw independent bootstrap innovations $\varepsilon_t^*$ from the empirical distribution $P_n$ given by the $\widehat{\varepsilon}_t$. A bootstrap version of the autoregressive process is now obtained as

$$X_t^* = \widehat{\theta}_{n,1} X_{t-1}^* + \cdots + \widehat{\theta}_{n,p} X_{t-p}^* + \varepsilon_t^*. \tag{2.5}$$

For simplicity, we assume that $(X_t^*)_{t\in\mathbb{Z}}$ is in its stationary regime. (This will be justified by (i) of the next lemma.) Before we state weak dependence of the bootstrap process, we show that it inherits those properties from the initial process which were used for proving weak dependence.

**Lemma 2.3.** *Suppose that $\mathbb{E}\varepsilon_0^2 < \infty$ and $\widehat{\theta}_n \xrightarrow{P} \theta$ are fulfilled.*

(i) *With a probability tending to 1, $(X_t^*)_{t\in\mathbb{Z}}$ can be written as a stationary causal linear process,*

$$X_t^* = \sum_{k=0}^{\infty} \widehat{\alpha}_{n,k} \varepsilon_{t-k}^*,$$

*where, for all $\epsilon > 0$, there exists a $\widetilde{K}_\epsilon < \infty$ such that, with $\rho_\epsilon = (1+\epsilon)/\rho$,*

$$\mathbb{P}\left(|\widehat{\alpha}_{n,k}| \leq \widetilde{K}_\epsilon \rho_\epsilon^k \quad \forall k \in \mathbb{N}\right) \xrightarrow[n\to\infty]{} 1.$$

(ii) $\mathbb{E}(\varepsilon_t^{*2} \mid X_{1-p}, \ldots, X_n) \xrightarrow{P} \mathbb{E}\varepsilon_t^2.$

Armed with the basic properties stated in Lemma 2.3, we can now easily derive properties of weak dependence of the bootstrap process just by imitating the proof for the initial process. In complete analogy to Lemma 2.2 above, we can state the following result.

**Lemma 2.4.** *Suppose that the initial process $(X_t)_{t\in\mathbb{Z}}$ satisfies the above conditions and that the bootstrap process $(X_t^*)_{t\in\mathbb{Z}}$ is in its stationary regime. Let $(X_t^{*'})_{t\in\mathbb{Z}}$ be another version of the bootstrap process, where $\mathbb{X}_{s_k}^* = (X_{s_k}^*, X_{s_k-1}^*, \ldots)$ is (conditionally on $X_{1-p}, \ldots, X_n$) independent of $\mathbb{X}_{s_k}^{*'} = (X_{s_k}^{*'}, X_{s_k-1}^{*'}, \ldots)$. For any $\epsilon > 0$, let $\rho_\epsilon = (1+\epsilon)/\rho$ and $\widetilde{\widetilde{K}}_\epsilon < \infty$ be an appropriate constant. Then there exists a sequence of events $\Omega_n$ such that $\mathbb{P}(\Omega_n) \xrightarrow[n\to\infty]{} 1$ and if $\Omega_n$ occurs, then*

$$\mathbb{E}\left(d\left(\mathbb{P}^{X_{t_1}^*,\ldots,X_{t_l}^* \mid \mathbb{X}_{s_k}^*}, \mathbb{P}^{X_{t_1}^{*'},\ldots,X_{t_l}^{*'} \mid \mathbb{X}_{s_k}^{*'}}\right) \mid X_{1-p}, \ldots, X_n\right) \leq 2\, l\, \widetilde{\widetilde{K}}_\epsilon\, \rho_\epsilon^{t_1-s_k}\, \sqrt{\mathbb{E}\varepsilon_0^2}.$$

From Lemma 2.1 and Lemma 2.4 we can now derive the desired property of $\psi$-weak dependence for the bootstrap process.

**Corollary 2.2.** *Suppose that the conditions of Lemma 2.4 are fulfilled. If the event $\Omega_n$ occurs, then the following assertion is true:*
*Let $s_1 \leq \cdots \leq s_k \leq t_1 \leq \cdots \leq t_l$ be arbitrary and let $g : \mathbb{R}^k \longrightarrow \mathbb{R}$ and $h :$*

$\mathbb{R}^l \longrightarrow \mathbb{R}$ *be measurable functions. If g is bounded and* $\mathbb{E}|h(X^*_{t_1},\ldots,X^*_{t_l})| < \infty$, *then*

$$\left|\operatorname{cov}\left(g(X^*_{s_1},\ldots,X^*_{s_k}), h(X^*_{t_1},\ldots,X^*_{t_l})\right)\right| \leq \|g\|_\infty \operatorname{Lip} h\; K_2\; l\; \rho_\epsilon^{t_1-s_k},$$

*where* $K_2 = 2\widetilde{\widetilde{K}}_\epsilon \sqrt{\mathbb{E}\varepsilon_0^2}$.

Besides the useful property of weak dependence of the bootstrap process, asymptotic validity of a bootstrap approximation requires that the (multivariate) stationary distributions of the bootstrap process converge to those of the initial process. Often, and in the case of the autoregressive bootstrap in particular, one has no direct access to these stationary distributions. However, according to Lemma 4.2 in Neumann and Paparoditis (2007), convergence of the stationary distributions can be derived from an appropriate convergence of conditional distributions. The latter, however, follows directly from $\widehat{\theta}_n \xrightarrow{P} \theta$ and $\widehat{\varepsilon}^*_t \xrightarrow{d} \varepsilon_t$. Therefore, consistency of the autoregressive bootstrap can be shown by simple arguments which were already used for proving weak dependence of the bootstrap; for details see Section 4.2 in Neumann and Paparoditis (2007).

*Remark* 3. Motivated by the desire to have some sort of mixing for a smoothed sieve bootstrap for linear processes, Bickel and Bühlmann (1999) considered a condition called $\nu$-mixing which is similar to the notion of weak dependence in our Definition 1.1. Although strong mixing follows for linear processes from a result of Gorodetskii (1977), it seems to be unclear whether even a smoothed version of the bootstrap process has such a property. However, it was shown in Theorems 3.2 and 3.4 in Bickel and Bühlmann (1999) that it is $\nu$-mixing with polynomial or exponential bounds on the corresponding coefficients to hold in probability. In the proofs of these theorems, however, they make use of the property of decaying strong mixing coefficients which holds at least for sufficiently large time lags; see in particular their Lemma 5.3.

In contrast, the approach described here is fundamentally different. We intend to prove weak dependence for processes driven by innovations with a possibly discrete distribution and achieve this goal by exploiting a contraction property of the initial and the bootstrap process.

*Remark* 4. Arguing in the same way as above we could also establish the property of $\psi$-weak dependence for nonlinear autoregressive processes,

$$X_t = m(X_t) + \varepsilon_t, \qquad t \in \mathbb{Z},$$

where $(\varepsilon_t)_{t\in\mathbb{Z}}$ is a sequence of independent and identically distributed innovations. If $\operatorname{Lip} m < 1$, then we have obviously a contraction property being fulfilled which immediately yields $\psi$-weak dependence.

It is interesting to note that such a contraction property can still be proved if $\operatorname{Lip} m < 1$ is not fulfilled. To this end, define the *local* Lipschitz modulus of continuity

$$\Delta(x) = \sup_{y\neq x} \frac{|m(y) - m(x)|}{|y - x|}$$

and assume that
$$\rho = \sup_x \mathbb{E}\left[\Delta(m(x) + \varepsilon_0)\right] < 1.$$

Then
$$d\left(\mathbb{P}^{X_{t+k}|X_t=x}, \mathbb{P}^{X_{t+k}|X_t=y}\right) \leq \rho^{k-1} \cdot \Delta(x) \cdot |x - y|, \tag{2.6}$$

which implies weak dependence by Lemma 2.1.

## 3. Some examples of weakly dependent sequences

Note first that sums of independent weakly dependent processes admit the common weak dependence property where dependence coefficients are the sums of the initial ones. We now provide a non-exhaustive list of weakly dependent sequences with their weak dependence properties. Further examples may be found in Doukhan and Louhichi (1999). Let $X = (X_t)_{t\in\mathbb{Z}}$ be a stationary process.

1. If this process is either a Gaussian process or an associated process and $\lim_{t\to\infty} |\text{cov}(X_0, X_t)| = 0$, then it is a $\kappa$-weakly dependent process such that $\kappa(r) = \mathcal{O}\left(\sup_{t\geq r} |\text{cov}(X_0, X_t)|\right)$. It is also $\kappa'$-weakly dependent with $\kappa'(r) = \mathcal{O}\left(\sum_{t\geq r} |\text{cov}(X_0, X_t)|\right)$.

2. $ARMA(p, q)$ processes and more generally causal or non-causal linear processes: $X = (X_t)_{t\in\mathbb{Z}}$ are defined by the model equation
$$X_t = \sum_{k=-\infty}^{\infty} a_k \xi_{t-k} \qquad \text{for } t \in \mathbb{Z},$$
where $(a_k)_{k\in\mathbb{Z}} \in \mathbb{R}^{\mathbb{Z}}$ and $(\xi_t)_{t\in\mathbb{Z}}$ is a sequence of independent and identically distributed random variables with $\mathbb{E}\xi_t = 0$. If $a_k = \mathcal{O}(|k|^{-\mu})$ with $\mu > 1/2$, then $X$ is an $\eta$-weakly dependent process with $\eta(r) = \mathcal{O}\left(\frac{1}{r^{\mu-1/2}}\right)$. In the general case of dependent innovations, properties of weak dependence are proved in Doukhan and Wintenberger (2007).

3. $GARCH(p, q)$ processes and more generally $ARCH(\infty)$ processes: $X = (X_t)_{t\in\mathbb{Z}}$ is a such that
$$X_t = \rho_t \cdot \xi_t \quad \text{with} \quad \rho_t^2 = b_0 + \sum_{k=1}^{\infty} b_k X_{t-k}^2 \quad \text{for } k \in \mathbb{Z},$$
with a sequence $(b_k)_k$ depending on the initial parameters in the case of a $GARCH(p, q)$ process and a sequence $(\xi_t)_{t\in\mathbb{Z}}$ of independent and identically distributed innovations. Then, if $\mathbb{E}(|\xi_0|^m) < \infty$, with the condition of stationarity, $\|\xi_0\|_m^2 \cdot \sum_{j=1}^{\infty} |b_j| < 1$, and if:

- there exists $C > 0$ and $\mu \in \,]0,1[$ such that $\forall j \in \mathbb{N}$, $0 \leq b_j \leq C \cdot \mu^{-j}$, then $X$ is a $\theta$-weakly dependent process with $\theta(r) = \mathcal{O}(e^{-c\sqrt{r}})$ and $c > 0$ (this is the case of $GARCH(p,q)$ processes).
- there exists $C > 0$ and $\nu > 1$ such that $\forall k \in \mathbb{N}$, $0 \leq b_k \leq C \cdot k^{-\nu}$, then $X$ is a $\theta$-weakly dependent process with $\theta(r) = \mathcal{O}(r^{-\nu+1})$ (Doukhan, Teyssière and Winant (2006) introduce vector valued LARCH($\infty$) models including the previous ones).

4. Causal bilinear processes were introduced by Giraitis and Surgailis (2002) and their dependence properties are proved in Doukhan, Madre and Rosenbaum (2007):

$$X_t = \xi_t \left( a_0 + \sum_{k=1}^{\infty} a_k X_{t-k} \right) + c_0 + \sum_{k=1}^{\infty} c_k X_{t-k}, \qquad \text{for } k \in \mathbb{Z}.$$

Assume that there exists $m \geq 1$ with $\|\xi_0\|_m < \infty$ and $\|\xi_0\|_m \cdot \left( \sum_{k=1}^{\infty} |a_k| + \sum_{k=1}^{\infty} |c_k| \right) < 1$. Then, if:

- $\begin{cases} \exists K \in \mathbb{N} \text{ such that } \forall k > K, a_k = c_k = 0, \text{ or,} \\ \exists \mu \in\, ]0,1[ \text{ such that } \sum_k |c_k|\mu^{-k} \leq 1 \text{ and } \forall k \in \mathbb{N}, 0 \leq a_k \leq \mu^k \end{cases}$
  then $X$ is a $\theta$-weakly dependent process with $\theta(r) = \mathcal{O}(e^{-c\sqrt{r}})$, for some $c > 0$;

- $\forall k \in \mathbb{N}$, $c_k \geq 0$, and $\exists \nu_1 > 2$ and $\exists \nu_2 > 0$ such that $a_k = \mathcal{O}(k^{-\nu_1})$ and $\sum_k c_k k^{1+\nu_2} < \infty$, then $X$ is a $\theta$-weakly dependent process with $\theta(r) = \mathcal{O}\left(\left(\frac{r}{\log r}\right)^d\right)$, $d = \max\left\{-(\nu_1 - 1); -\frac{\nu_2 \cdot \delta}{\delta + \nu_2 \cdot \log 2}\right\}$ (see Doukhan, Teyssière and Winant (2006)).

5. Non-causal LARCH($\infty$) processes $X = (X_t)_{t \in \mathbb{Z}}$ satisfying

$$X_t = \xi_t \cdot \left( a_0 + \sum_{k \in \mathbb{Z} \setminus \{0\}} a_k X_{t-k} \right), \qquad t \in \mathbb{Z},$$

where $\|\xi_0\|_\infty < \infty$ (bounded random variables) and $(a_k)_{k \in \mathbb{Z}}$ is a sequence of real numbers such that $\lambda = \|\xi_0\|_\infty \cdot \sum_{j \neq 0} |a_j| < 1$ (stationarity condition). Assume that the sequence $(a_k)_{k \in \mathbb{Z}}$ satisfies $a_k = \mathcal{O}(|k|^{-\mu})$ with $\mu > 1$, then $X$ is an $\eta$-weakly dependent process with $\eta(r) = O\left(\frac{1}{r^{\mu-1}}\right)$ (see Doukhan, Teyssière and Winant (2006)).

6. Causal and non-causal Volterra processes write as $X_t = \sum_{p=1}^{\infty} Y_t^{(p)}$ with

$$Y_t^{(p)} = \sum_{\substack{j_1 < j_2 < \cdots < j_p \\ j_1, \ldots, j_p \in \mathbb{Z}}} a_{j_1, \ldots, j_p} \xi_{t-j_1} \cdots \xi_{t-j_p}, \qquad \text{for } t \in \mathbb{Z}.$$

Assume $\sum_{p=0}^{\infty} \sum_{\substack{j_1 < j_2 < \cdots < j_p \\ j_1, \ldots, j_p \in \mathbb{Z}}} |a_{j_1, \ldots, j_p}|^m \|\xi_0\|_m^p < \infty$, with $m > 0$, and that there exists $p_0 \in \mathbb{N} \setminus \{0\}$ such that $a_{j_1, \ldots, j_p} = 0$ for $p > p_0$. If $a_{j_1, \ldots, j_p} = \mathcal{O}\left(\max_{1 \leq i \leq p}\{|j_i|^{-\mu}\}\right)$ with $\mu > 0$, then $X$ is an $\eta$-weakly dependent process with $\eta(r) = \mathcal{O}\left(\frac{1}{r^{\mu+1}}\right)$ (see Doukhan (2002)).

Finite order Volterra processes with dependent inputs are also considered in Doukhan and Wintenberger (2007): again, $\eta$-weakly dependent innovations yield $\eta$-weak dependence and $\lambda$-weakly dependent innovation yields $\lambda$-weak dependence of the process.

7. Very general models are the causal or non-causal infinite memory processes $X = (X_t)_{t \in \mathbb{Z}}$ such that

$$X_t = F(X_{t-1}, X_{t-2}, \ldots; \xi_t), \quad \text{and} \quad X_t = F(X_s, s \neq t; \xi_t),$$

where the functions $F$ defined either on $\mathbb{R}^{\mathbb{N} \setminus \{0\}} \times \mathbb{R}$ or $\mathbb{R}^{\mathbb{Z} \setminus \{0\}} \times \mathbb{R}$ satisfy

$$\|F(0; \xi_0)\|_m < \infty,$$
$$\|F((x_j)_j; \xi_0) - F((y_j)_j; \xi_0)\|_m \leq \sum_{j \neq 0} a_j |x_j - y_j|,$$

with $a = \sum_{j \neq 0} a_j < 1$. Then, works in progress by Doukhan and Wintenberger as well as Doukhan and Truquet, respectively, prove that a solution of the previous equations is stationary in $L^m$ and either $\theta$-weakly dependent or $\eta$-weakly dependent with the following decay rate for the coefficients:

$$\inf_{p \geq 1} \left\{ a^{r/p} + \sum_{|j| > p} a_j \right\}.$$

This provides the same rates as those already mentioned for the cases of $ARCH(\infty)$ or $LARCH(\infty)$ models.

## 4. Some probabilistic results

In this section, we present results derived under weak dependence which are of interest in probability and statistics (see also Dedecker *et al.* (2007) for reference). This collection clearly shows that this notion of weak dependence, although being more general than mixing, allows one to prove results very similar to those in the mixing case.

### *4.1. Donsker invariance principle*

We consider a stationary, zero mean, and real valued sequence $(X_t)_{t \in \mathbb{Z}}$ such that

$$\mu = \mathbb{E}|X_0|^m < \infty, \quad \text{for a real number } m > 2. \quad (4.1)$$

We also set

$$\sigma^2 = \sum_{t\in\mathbb{Z}} \operatorname{cov}(X_0, X_t) = \sum_{t\in\mathbb{Z}} \mathbb{E} X_0 X_t, \tag{4.2}$$

$W$ denotes standard Brownian motion and

$$W_n(t) = \frac{1}{\sqrt{n}} \sum_{i=1}^{[nt]} X_i, \quad \text{for } t \in [0,1],\ n \geq 1. \tag{4.3}$$

We now present versions of the Donsker weak invariance principle under weak dependence assumptions.

**Theorem 4.1** (Donsker type results). *Assume that the zero mean stationary process $(X_t)_{t\in\mathbb{Z}}$ satisfies (4.1). Then $\sigma^2 \geq 0$ given by (4.2) is well defined and the Donsker invariance principle*

$W_n(t) \to_{n\to\infty} \sigma W(t), \quad$ *in distribution in the Skorohod space $D([0,1])$,*

*holds if one of the following additional assumptions is fulfilled:*

- **$\kappa$-dependence.** *The process is $\kappa$-weakly dependent and satisfies $\kappa(r) = \mathcal{O}(r^{-\kappa})$ (as $r \uparrow \infty$) for some $\kappa > 2 + 1/(m-2)$.*
- **$\kappa'$-dependence.** *The process is $\kappa'$-weakly dependent and satisfies $\kappa'(r) = \mathcal{O}(r^{-\kappa})$ (as $r \uparrow \infty$) for some $\kappa > 1 + 1/(m-2)$.*
- **$\lambda$-dependence.** *The process is $\lambda$-weakly dependent and satisfies $\lambda(r) = \mathcal{O}(r^{-\lambda})$ (as $r \uparrow \infty$) for $\lambda > 4 + 2/(m-2)$.*
- **$\theta$-dependence.** *The process is $\theta$-weakly dependent and satisfies $\theta(r) = \mathcal{O}(r^{-\theta})$ (as $r \uparrow \infty$) for $\theta > 1 + 1/(m-2)$.*

*Remark* 5. The result for $\kappa'$-weak dependence is obtained in Bulinski and Shashkin (2005). Results under $\kappa$- and $\lambda$-weak dependences are proved in Doukhan and Wintenberger (2007); note that $\eta$-weak dependence implies $\lambda$-weak dependence and the Donsker principle then holds under the same decay rate for the coefficients. The result for $\theta$-weak dependence is due to Dedecker and Doukhan (2003). A few comments on these results are now in order:

- The difference of the above conditions under $\kappa$ and $\kappa'$ assumptions is natural. The observed loss under $\kappa$-dependence is explained by the fact that $\kappa'$-weakly dependent sequences satisfy $\kappa'(r) \geq \sum_{s \geq r} \kappa(s)$. This simple bound directly follows from the definitions.
- Actually, it is enough to assume the $\theta$-weak dependence inequality for any positive integer $u$ and only for $v = 1$. Hence, for any 1-bounded function $g$ from $\mathbb{R}^u$ to $\mathbb{R}$ and any 1-bounded Lipschitz function $h$ from $\mathbb{R}$ to $\mathbb{R}$ with Lipschitz coefficient $\operatorname{Lip}(h)$, it is enough to assume that the following inequality is fulfilled: $\left|\operatorname{cov}\left(g(X_{i_1},\ldots,X_{i_u}), h(X_{i_u+i})\right)\right| \leq \theta(i)\operatorname{Lip}(h)$, for any $u$-tuple $i_1 \leq i_2 \leq \cdots \leq i_u$.

### 4.2. Empirical process

Let $(X_t)_{t\in\mathbb{Z}}$ a real-valued stationary process. We use a quantile transform to obtain that the marginal distribution of this sequence is the uniform law on $[0,1]$. The empirical process of the sequence $(X_t)_{t\in\mathbb{Z}}$ at time $n$ is defined as $\frac{1}{\sqrt{n}}E_n(x)$ where

$$E_n(x) = \sum_{k=1}^{n} \left(\mathbb{1}_{(X_k \leq x)} - \mathbb{P}(X_k \leq x)\right).$$

Note that $E_n = n(F_n - F)$ if $F_n$ and $F$ denote the empirical distribution function and the marginal distribution function, respectively. We consider the following convergence result in the Skohorod space $D([0,1])$ when the sample size $n$ tends to infinity:

$$\frac{1}{\sqrt{n}}E_n(x) \xrightarrow{d} \bar{B}(x).$$

Here $(\bar{B}(x))_{x\in[0,1]}$ is the dependent analogue of a Brownian bridge, that is $\bar{B}$ denotes a centered Gaussian process with covariance given by

$$\mathbb{E}\bar{B}(x)\bar{B}(y) = \sum_{k=-\infty}^{\infty} \left(\mathbb{P}(X_0 \leq x, X_k \leq y) - \mathbb{P}(X_0 \leq x)\mathbb{P}(X_k \leq y)\right). \quad (4.4)$$

Note that for independent sequences with a marginal distribution function $F$, this turns into $\bar{B}(x) = B(x)$ for some standard Brownian bridge $B$; this justifies the name of generalized Brownian bridge. We have:

**Theorem 4.2.** *Suppose that the stationary sequence $(X_t)_{t\in\mathbb{Z}}$ has a uniform marginal distribution and is $\eta$-weakly dependent with $\eta(r) = \mathcal{O}(r^{-15/2-\nu})$, or $\kappa$-weakly dependent $\kappa(r) = \mathcal{O}(r^{-5-\nu})$, for some $\nu > 0$. Then the following empirical functional convergence holds true in the Skohorod space of real-valued càdlàg functions on the unit interval, $D([0,1])$:*

$$\frac{1}{\sqrt{n}}E_n(x) \xrightarrow{d} \bar{B}(x).$$

*Remark* 6. Under strong mixing, the condition $\sum_{r=0}^{\infty} \alpha(r) < \infty$ implies convergence of the finite-dimensional distributions. The empirical functional convergence holds if, in addition, for some $a > 1$, $\alpha(r) = \mathcal{O}(r^{-a})$ (see Rio (2000)). In an absolutely regular framework, Doukhan, Massart and Rio (1995) obtain the empirical functional convergence when, for some $a > 2$, $\beta(r) = \mathcal{O}(r^{-1}(\log r)^{-a})$. Shao and Yu (1996) and Shao (1995) obtain the empirical functional convergence theorem when the maximal correlation coefficients satisfy the condition $\sum_{n=0}^{\infty} \rho(2^n) < \infty$.

To prove the result, we introduce the following dependence condition for a stationary sequence $(X_t)_{t\in\mathbb{Z}}$:

$$\sup_{f\in\mathcal{F}} |\mathrm{cov}(f(X_{t_1})f(X_{t_2}), f(X_{t_3})f(X_{t_4}))| \leq \epsilon(r), \quad (4.5)$$

where $\mathcal{F} = \{x \mapsto \mathbb{I}_{s<x\leq t}, \text{ for } s, t \in [0,1]\}$, $0 \leq t_1 \leq t_2 \leq t_3 \leq t_4$ and $r = t_3 - t_2$ (in this case a weak dependence condition holds for a class of functions $\mathbb{R}^u \to \mathbb{R}$ working only with the values $u = 1$ or $2$).

**Proposition 4.1.** *Let $(X_n)$ be a stationary sequence such that (4.5) holds. Assume that there exists $\nu > 0$ such that*

$$\epsilon(r) = \mathcal{O}(r^{-5/2-\nu}). \tag{4.6}$$

*Then the sequence of processes $\left(\left\{\frac{1}{\sqrt{n}}E_n(t); t \in [0,1]\right\}\right)_{n>0}$ is tight in the Skohorod space $D([0,1])$.*

### 4.3. Central limit theorems

First central limit theorems for weakly dependent sequences were given by Corollary A in Doukhan and Louhichi (1999) and Theorem 1 in Coulon-Prieur and Doukhan (2000). While the former result is for sequences of stationary random variables, the latter one is tailor-made for triangular arrays of asymptotically sparse random variables as they appear with kernel density estimators. Using their notion of $\nu$-mixing Bickel and Bühlmann (1999) proved a CLT for linear processes of infinite order and their (smoothed) bootstrap counterparts. Below we state a central limit theorem for general triangular schemes of weakly dependent random variables. Note that the applicability of a central limit theorem to bootstrap processes requires some robustness in the parameters of the underlying process since these parameters have to be estimated when it comes to the bootstrap. A result for a triangular scheme is therefore appropriate since the involved random variables have themselves random properties. An interesting aspect of the following results is that no moment condition beyond Lindeberg's is required. The proof of the next theorem uses the variant of Rio of the classical Lindeberg method.

**Theorem 4.3.** *(Theorem 6.1 in Neumann and Paparoditis (2007)) Suppose that $(X_{n,k})_{k=1,\ldots,n}$, $n \in \mathbb{N}$, is a triangular scheme of (row-wise) stationary random variables with $\mathbb{E}X_{n,k} = 0$ and $\mathbb{E}X_{n,k}^2 \leq C < \infty$. Furthermore, we assume that*

$$\frac{1}{n}\sum_{k=1}^{n}\mathbb{E}X_{n,k}^2 I(|X_{n,k}|/\sqrt{n} > \epsilon) \underset{n\to\infty}{\longrightarrow} 0 \tag{4.7}$$

*holds for all $\epsilon > 0$ and that*

$$\mathrm{var}(X_{n,1} + \cdots + X_{n,n})/n \underset{n\to\infty}{\longrightarrow} \sigma^2 \in [0,\infty). \tag{4.8}$$

*For $n \geq n_0$, there exists a monotonously nonincreasing and summable sequence $(\theta(r))_{r\in\mathbb{N}}$ such that, for all indices $s_1 < s_2 < \cdots < s_u < s_u + r = t_1 \leq t_2$, the following upper bounds for covariances hold true: for all measurable and quadratic integrable functions $g : \mathbb{R}^u \to \mathbb{R}$,*

$$|\mathrm{cov}\,(g(X_{n,s_1},\ldots,X_{n,s_u}), X_{n,t_1})| \leq \sqrt{\mathbb{E}g^2(X_{n,s_1},\ldots,X_{n,s_u})}\,\theta(r), \tag{4.9}$$

*for all measurable and bounded functions* $g : \mathbb{R}^u \longrightarrow \mathbb{R}$,

$$|\text{cov}(g(X_{n,s_1}, \ldots, X_{n,s_u}), X_{n,t_1} X_{n,t_2})| \leq \|g\|_\infty \, \theta(r), \quad (4.10)$$

*where* $\|g\|_\infty = \sup_{x \in \mathbb{R}^u} |f(x)|$. *Then*

$$\frac{1}{\sqrt{n}}(X_{n,1} + \cdots + X_{n,n}) \xrightarrow{d} \mathcal{N}(0, \sigma^2).$$

*Remark 7.* (i) Conditions (4.7) to (4.10) can be easily verified for autoregressive processes under the standard condition that the characteristic polynomial has no zero within the unit circle and if $\mathbb{E}\varepsilon_0^2 < \infty$, that is, finiteness of second moments actually suffices here; see Neumann and Paparoditis (2007) for details.

(ii) If we have in the $n$th stage a two-sided sequence of random variables $(X_{n,k})_{k \in \mathbb{Z}}$ rather than $(X_{n,k})_{k=1,\ldots,n}$ only, then it can be easily seen that condition (4.8) follows from $EX_{n,k} = 0$, $EX_{n,k}^2 \leq C < \infty$, (4.9), and $\sum_{k \in \mathbb{Z}} \text{cov}(X_{n,0}, X_{n,k}) \xrightarrow[n \to \infty]{} \sigma^2 \in [0, \infty)$.

(iii) Condition (4.9) is also related to Gordin (1969)'s condition under which central limit theorems are often proved for stationary processes. Such a theorem for a sequence of stationary ergodic random variables was proved by Hall and Heyde (1980, pp. 136–138); see also Esseen and Janson (1985) for the correction of a detail.

The following very simple multivariate central limit theorem, easily applicable to triangular schemes of weakly dependent random vectors, was derived in Bardet, Doukhan, Lang and Ragache (2007). In view of condition (4.11), it is applicable in cases where dependence between the observations declines as $n \to \infty$. This is a common situation in nonparametric curve estimation where the so-called "whitening-by-windowing" principle applies.

**Theorem 4.4.** *(Theorem 1 in Bardet, Doukhan, Lang and Ragache (2007)) Suppose that $(X_{n,k})_{k \in \mathbb{N}}$, $n \in \mathbb{N}$, is a triangular scheme of zero mean random vectors with values in $\mathbb{R}^d$. Assume that there exists a positive definite matrix $\Sigma$ such that*

$$\sum_{k=1}^{n} \text{Cov}(X_{n,k}) \xrightarrow[n \to \infty]{} \Sigma$$

*and that, for each $\epsilon > 0$,*

$$\sum_{k=1}^{n} \mathbb{E}(\|X_{n,k}\|^2 \mathbb{1}_{\{\|X_{n,k}\| > \epsilon\}}) \xrightarrow[n \to \infty]{} 0,$$

*where $\|\cdot\|$ denotes the Euclidean norm. Furthermore, we assume the following condition is satisfied:*

$$\sum_{k=2}^{n} \left|\text{cov}(e^{it'(X_{n,1} + \cdots + X_{n,k-1})}, e^{it' X_{n,k}})\right| \xrightarrow[n \to \infty]{} 0. \quad (4.11)$$

*Then, as $n \to \infty$,*

$$\sum_{k=1}^{n} X_{n,k} \xrightarrow{d} \mathcal{N}_d(0_d, \Sigma).$$

*Remark* 8. One common point of these two results is the use of the classical Lindeberg assumption. Note that this assumption is often verified by using a higher order moment condition. A main difference between the two results is that the first one yields direct applications to partial sums while the second one is more adapted to triangular arrays where the limit does not write as a sum. In this setting Doukhan and Wintenberger (2007) use Bernstein blocks to prove a CLT for partial sums.

### 4.4. Probability and moment inequalities

In this section we state inequalities of Bernstein and Rosenthal type. In the case of mixing, such inequalities can be easily derived by the well-known technique of replacing dependent blocks of random variables (separated by an appropriate time gap) by independent ones and then using the classical inequalities from the independent case; see for example Doukhan (1994) and Rio (2000). The notion of $\psi$-weak dependence is particularly suitable for deriving upper estimates for the cumulants of sums of random variables which give rise to rather sharp inequalities of Bernstein and Rosenthal type which are analogous to those in the independent case.

Based on a Rosenthal-type inequality, a first inequality of Bernstein-type was obtained by Doukhan and Louhichi (1999), however, with $\sqrt{t}$ instead of $t^2$ in the exponent. Dedecker and Prieur (2004) proved a Bennett inequality which can possibly be used to derive also a Bernstein inequality. A first Bernstein inequality with $\text{var}(X_1 + \cdots + X_n)$ in the asymptotically leading term of the denominator of the exponent has been derived in Kallabis and Neumann (2006), under a weak dependence condition tailor-made for causal processes with an exponential decay of the coefficients of weak dependence. The following result is a generalization which is also applicable to possibly non-causal processes with a not necessarily exponential decay of the coefficients of weak dependence.

**Theorem 4.5.** *(Theorem 1 in Doukhan and Neumann (2007)) Suppose that $X_1, \ldots, X_n$ are real-valued random variables with zero mean, defined on a probability space $(\Omega, \mathcal{A}, \mathbb{P})$. Let $\Psi : \mathbb{N}^2 \to \mathbb{N}$ be one of the following functions:*

*(a)* $\Psi(u, v) = 2v$,
*(b)* $\Psi(u, v) = u + v$,
*(c)* $\Psi(u, v) = uv$,
*(d)* $\Psi(u, v) = \alpha(u + v) + (1 - \alpha)uv$, *for some $\alpha \in (0, 1)$.*

*We assume that there exist constants $K, M, L_1, L_2 < \infty$, $\mu, \nu \geq 0$, and a nonincreasing sequence of real coefficients $(\rho(n))_{n \geq 0}$ such that, for all $u$-tuples* $(s_1, \ldots, s_u)$ *and all v-tuples* $(t_1, \ldots, t_v)$ *with* $1 \leq s_1 \leq \cdots \leq s_u \leq t_1 \leq \cdots \leq t_v \leq n$ *the following inequalities are fulfilled:*

$$|\text{cov}(X_{s_1} \cdots X_{s_u}, X_{t_1} \cdots X_{t_v})| \leq K^2 \, M^{u+v-2} \, ((u+v)!)^\nu \, \Psi(u,v) \, \rho(t_1 - s_u), \quad (4.12)$$

*where*

$$\sum_{s=0}^{\infty} (s+1)^k \rho(s) \leq L_1 \, L_2^k \, (k!)^\mu \qquad \forall k \geq 0, \quad (4.13)$$

*and*

$$\mathbb{E}|X_t|^k \leq (k!)^\nu \, M^k \qquad \forall k \geq 0. \quad (4.14)$$

*Then, for all* $t \geq 0$,

$$\mathbb{P}(S_n \geq t) \leq \exp\left(-\frac{t^2/2}{A_n + B_n^{\frac{1}{\mu+\nu+2}} t^{\frac{2\mu+2\nu+3}{\mu+\nu+2}}}\right), \quad (4.15)$$

*where* $A_n$ *can be chosen as any number greater than or equal to* $\sigma_n^2$ *and*

$$B_n = 2 \, (K \vee M) \, L_2 \left(\left(\frac{2^{4+\mu+\nu} \, n \, K^2 \, L_1}{A_n}\right) \vee 1\right).$$

A first Rosenthal-type inequality for weakly dependent random variables was derived by Doukhan and Louhichi (1999) via direct expansions of the moments of even order. Unfortunately, the variance of the sum did not explicitly show up in their bound. Using cumulant bounds in conjunction with Leonov and Shiryaev's formula the following tighter moment inequality was obtained in Doukhan and Neumann (2007).

**Theorem 4.6.** *(Theorem 3 in Doukhan and Neumann (2007)) Suppose that* $X_1, \ldots, X_n$ *are real-valued random variables on a probability space* $(\Omega, \mathcal{A}, \mathbb{P})$ *with zero mean and let p be a positive integer. We assume that there exist constants* $K, M < \infty$, *and a non-increasing sequence of real coefficients* $(\rho(n))_{n \geq 0}$ *such that, for all u-tuples* $(s_1, \ldots, s_u)$ *and all v-tuples* $(t_1, \ldots, t_v)$ *with* $1 \leq s_1 \leq \cdots \leq s_u \leq t_1 \leq \cdots \leq t_v \leq n$ *and* $u + v \leq p$, *condition (4.12) is fulfilled. Furthermore, we assume that*

$$\mathbb{E}|X_i|^{p-2} \leq M^{p-2}.$$

*Then, with* $Z \sim \mathcal{N}(0,1)$,

$$\left|\mathbb{E}\left(\sum_{k=1}^n X_k\right)^p - \sigma_n^p \mathbb{E} Z^p\right| \leq B_{p,n} \sum_{1 \leq u < p/2} A_{u,p} \, K^{2u} \, (M \vee K)^{p-2u} \, n^u,$$

*where* $B_{p,n} = (p!)^2 2^p \max_{2 \leq k \leq p} \{\rho_{k,n}^{p/k}\}$, $\rho_{k,n} = \sum_{s=0}^{n-1} (s+1)^{k-2} \rho(s)$ *and*

$$A_{u,p} = \frac{1}{u!} \sum_{k_1 + \cdots + k_u = p, \, k_i \geq 2 \, \forall i} \frac{p!}{k_1! \cdots k_u!}.$$

To avoid any misinterpretation, we note that condition (4.12) with $u + v \leq p$ and $\mathbb{E}|X_i|^{p-2} \leq M^{p-2}$ only requires finiteness of moments of order $p$. This is in contrast to the conditions imposed in Theorem 4.5 where in particular all moments of the involved random variables have to be finite.

## 5. Proofs

*Proof of Lemma 2.1.*
Let, for simplicity of notation, $s_k = 0$. Then

$$\text{cov}\,(g(X_{s_1},\ldots,X_{s_k}), h(X_{t_1},\ldots,X_{t_l}))$$
$$= \mathbb{E}\left[g(X_{s_1},\ldots,X_{s_k})\left(\mathbb{E}(h(X_{t_1},\ldots,X_{t_l}) \mid X_{s_1},\ldots,X_{s_k}) - \mathbb{E}h(X_{t_1},\ldots,X_{t_l})\right)\right].$$

Now we obtain by Jensen's inequality for conditional expectations that

$$|\text{cov}\,(g(X_{s_1},\ldots,X_{s_k}), h(X_{t_1},\ldots,X_{t_l}))|$$
$$\leq \mathbb{E}\left[|g(X_{s_1},\ldots,X_{s_k})| \cdot |\mathbb{E}(h(X_{t_1},\ldots,X_{t_l}) \mid \mathbb{X}_0) - \mathbb{E}h(X_{t_1},\ldots,X_{t_l})|\right]$$
$$\leq \mathbb{E}\left[|g(X_{s_1},\ldots,X_{s_k})| \cdot \left|\mathbb{E}(h(X_{t_1},\ldots,X_{t_l}) \mid \mathbb{X}_0) - \mathbb{E}(h(X'_{t_1},\ldots,X'_{t_l}) \mid \mathbb{X}'_0)\right|\right]$$
$$\leq \text{Lip}\,h\, \mathbb{E}\left[|g(X_{s_1},\ldots,X_{s_k})| \cdot \mathbb{E}\left(\sum_{j=1}^{l} |X_{t_j} - X'_{t_j}| \mid \mathbb{X}_0, \mathbb{X}'_0\right)\right].$$

The assertion follows now immediately. □

*Proof of Lemma 2.2.*
The assertion follows immediately from (2.3) and (2.4). □

*Proof of Lemma 2.3.*
(i) Let $\widehat{\xi}_{n,1},\ldots,\widehat{\xi}_{n,p}$ be the roots of the characteristic polynomial $\widehat{\theta}_n(z) = 1 - \widehat{\theta}_{n,1}z - \cdots - \widehat{\theta}_{n,p}z^p$ of the bootstrap process. Since $\widehat{\theta}_n \xrightarrow{P} \theta$ we obtain by Theorem 1.4 in Marden (1949) that

$$\min\{|\widehat{\xi}_{n,1}|,\ldots,|\widehat{\xi}_{n,p}|\} \xrightarrow{P} \rho = \min\{|\xi_1|,\ldots,|\xi_p|\}.$$

Therefore, we have, for any $\epsilon > 0$, that

$$\mathbb{P}\left(\min\{|\widehat{\xi}_{n,1}|,\ldots,|\widehat{\xi}_{n,p}|\} \geq \rho/(1 + \epsilon/2)\right) \xrightarrow[n\to\infty]{} 1. \qquad (5.1)$$

Thus there exists a stationary solution to equation (2.5) which can also be written as a causal linear process,

$$X_t^* = \sum_{k=0}^{\infty} \widehat{\alpha}_{n,k}\varepsilon_{t-k}^*.$$

It follows from (5.1) that, for all $\epsilon > 0$, there exist some $K_\epsilon < \infty$ such that, with $\rho_\epsilon = (1+\epsilon)/\rho$,

$$\mathbb{P}\left(|\widehat{\alpha}_{n,k}| \leq C_\epsilon \rho_\epsilon^k \quad \forall k\right) \xrightarrow[n\to\infty]{} 1.$$

(ii) To prove (ii), we first split up

$$\frac{1}{n}\sum_{t=1}^n \widetilde{\varepsilon}_t^2 = \frac{1}{n}\sum_{t=1}^n \varepsilon_t^2 + \frac{1}{n}\sum_{t=1}^n (\widetilde{\varepsilon}_t - \varepsilon_t)^2 + \frac{2}{n}\sum_{t=1}^n (\widetilde{\varepsilon}_t - \varepsilon_t)\varepsilon_t$$
$$= T_{n,1} + T_{n,2} + T_{n,3},$$

say. It follows from the strong law of large numbers that

$$T_{n,1} \xrightarrow{a.s.} \mathbb{E}\varepsilon_0^2.$$

Since $(X_t)_{t\in\mathbb{Z}}$ is ergodic we obtain from Birkhoff (1931)'s ergodic theorem (see also Corollary 3.5.1 in Stout (1974)) that $\frac{1}{n}\sum_{t=1}^n \boldsymbol{X}_t \boldsymbol{X}_t' \xrightarrow{a.s.} \mathbb{E}\boldsymbol{X}_0\boldsymbol{X}_0'$. Using $\widetilde{\varepsilon}_t - \varepsilon_t = \boldsymbol{X}_t'(\theta - \widehat{\theta}_n)$ we therefore obtain that

$$T_{n,2} = (\widehat{\theta}_n - \theta)'\frac{1}{n}\sum_{t=1}^n \boldsymbol{X}_t\boldsymbol{X}_t'(\widehat{\theta}_n - \theta) \xrightarrow{P} 0.$$

Finally, we conclude by the Cauchy-Schwarz inequality that $T_{n,3} \xrightarrow{P} 0$, which gives that

$$\frac{1}{n}\sum_{t=1}^n \widetilde{\varepsilon}_t^2 \xrightarrow{P} \mathbb{E}\varepsilon_0^2.$$

Since by the strong law of large numbers $\overline{\varepsilon}. = (1/n)\sum_{t=1}^n \varepsilon_t \xrightarrow{a.s.} 0$ and by Jensen's inequality $(\widetilde{\varepsilon}. - \overline{\varepsilon}.)^2 \leq (1/n)\sum_{t=1}^n (\widetilde{\varepsilon}_t - \varepsilon_t)^2$ we also obtain that $\widetilde{\varepsilon}. \xrightarrow{P} 0$. This implies that

$$\mathbb{E}\left(\varepsilon_t^* \mid X_{1-p},\ldots,X_n\right) = \frac{1}{n}\sum_{t=1}^n \widehat{\varepsilon}_t^2 = \frac{1}{n}\sum_{t=1}^n \widetilde{\varepsilon}_t^2 - \widetilde{\varepsilon}.^2 \xrightarrow{P} \mathbb{E}\varepsilon_0^2.$$

□

*Proof of inequality (2.6).*
We prove this result by a simple coupling argument. Let $(X_t)_{t\in\mathbb{Z}}$ and $(X_t')_{t\in\mathbb{Z}}$ be two versions of the autoregressive process with $X_t = x$ and $X_t' = y$. We contruct a coupling simply by feeding both processes after time $t$ with the same sequence of innovations $\varepsilon_{t+1}, \varepsilon_{t+2}, \ldots$, that is, we have $X_{t+l+1} = m(X_{t+l}) + \varepsilon_{t+l+1}$ and $X_{t+l+1}' = m(X_{t+l}') + \varepsilon_{t+l+1}$ ($l \geq 0$). It follows from this construction that

$$|X_{t+k} - X_{t+k}'| \leq \Delta(X_{t+k-1})|X_{t+k-1} - X_{t+k-1}'|$$
$$\leq \cdots$$
$$\leq \Delta(X_{t+k-1})\cdots\Delta(X_{t+1})\Delta(x)|x-y|.$$

Therefore, we obtain that

$$d\left(\mathbb{P}^{X_{t+k}|X_t=x}, \mathbb{P}^{X_{t+k}|X_t=y}\right)$$
$$\leq \mathbb{E}|X_{t+k} - X_{t+k}'|$$
$$\leq \mathbb{E}\left(\Delta(X_{t+k-1})\cdots\Delta(X_{t+1})\mid X_t = x\right)\cdot\Delta(x)\cdot|x-y|.$$

Using the Markov property we see that

$$\begin{aligned}
&\mathbb{E}\left(\Delta(X_{t+k-1})\cdots\Delta(X_{t+1})|\, X_t = x\right) \\
&= \mathbb{E}\left[\mathbb{E}(\Delta(X_{t+k-1}) \mid X_{t+k-2})\cdots\mathbb{E}(\Delta(X_{t+1}) \mid X_t = x)\right] \\
&\leq \rho^{k-1},
\end{aligned}$$

which yields the assertion. □

*Acknowledgment*. We thank an Associate Editor and an anonymous referee for their helpful comments and suggestions.